\documentclass[11pt]{article}
\usepackage{amsmath}
\usepackage{setspace}
\usepackage{fullpage}
\usepackage{graphicx}
\usepackage{amsthm}

\usepackage{epsfig}
\usepackage{epsf}  

\newtheorem{thm}{Theorem}[section]
\newtheorem{lemma}[thm]{Lemma}
\newtheorem{prop}[thm]{Proposition}
\newtheorem{corr}[thm]{Corollary}

\theoremstyle{definition}

\begin{document}
 
\title{On $3$-regular $4$-ordered graphs }
 
\author{
Karola M\'esz\'aros\\
Massachusetts  Institute of Technology\\
 {\tt karola@math.mit.edu}
\\
%\vspace{1in}
}
\date{}
\maketitle

\begin{abstract}
 
\begin{small} 
A simple graph $G$ is \textit{k-ordered} (respectively, \textit{k-ordered
hamiltonian}), if for any sequence of $k$ distinct vertices $v_1, \ldots, v_k$
of $G$ there exists a cycle (respectively, hamiltonian cycle) in $G$
containing these $k$ vertices in the specified order.  In 1997 Ng and Schultz
introduced  these concepts of cycle orderability and posed the
question of the existence of $3$-regular $4$-ordered (hamiltonian) graphs
other than $K_4$ and $K_{3, 3}$. Ng and Schultz observed that a $3$-regular
$4$-ordered graph on more than $4$ vertices is triangle free.  We prove that a
$3$-regular $4$-ordered graph $G$ on more than $6$ vertices is square free,
and we show that the smallest graph that is triangle and square free, namely
the Petersen graph, is $4$-ordered.  Furthermore, we prove that the smallest
graph after $K_4$ and $K_{3, 3}$ that is $3$-regular $4$-ordered hamiltonian
is the Heawood graph, and we exhibit forbidden subgraphs for $3$-regular
$4$-ordered hamiltonian graphs on more than $10$ vertices.  Finally, we
construct an infinite family of $3$-regular $4$-ordered graphs. 
 \end{small}
\end{abstract}
 
\section{Introduction}

The concept of $k$-ordered graphs was introduced in 1997 by Ng and
Schultz \cite{ng}. A simple graph  $G$ is a graph without loops or multiple
edges, and it is called  \textit{hamiltonian} if there exists a cycle that
contains all vertices of $G$.   In this paper we consider only connected
finite simple graphs.  A simple graph $G$ is called \textit{k-ordered}
(respectively, \textit{k-ordered hamiltonian}), if for any sequence of $k$
distinct vertices $v_1, \ldots, v_k$ of $G$ there exists a cycle
(respectively, hamiltonian cycle) in $G$ containing these $k$ vertices in the
specified order. Previous results concerning cycle orderability focus on
minimum degree and forbidden subgraph conditions that imply $k$-orderedness or
$k$-ordered hamiltonicity \cite{chen, faun-faun, faun-gould}.  A comprehensive
survey of results can be found in \cite{faun}.

Any hamiltonian graph is necessarily $3$-ordered hamiltonian as well as
$3$-ordered, thus we study $k$-orderedness for $k\geq4$. Indeed, it is easy
to see that hamiltonicity is equivalent to $3$-ordered hamiltonicity, and
$3$-cyclability to $3$-orderedness  (a graph is said to be
$3$-\textit{cyclable}, if for any three vertices of the graph there exists a cycle
containing them).  If $G$ is a  hamiltonian graph of order $n\geq 3$, then $G$ being
$k$-ordered
hamiltonian implies that $G$ is $(k-1)$-connected (see \cite{ng}). 
The arguments made in \cite{ng} hold in case of $k$-orderedness as well,  
namely, if $G$ is a graph of order $n\geq 3$, then $G$ being
$k$-ordered implies that $G$ is $(k-1)$-connected.  In particular, this
implies that $\delta(G)$, the minimum  degree of any vertex in a $k$-ordered graph $G$ 
is at least $k-1$. 
 
In \cite{ng}, a search for small degree $4$-ordered hamiltonian graphs
was conducted and an infinite family of $4$-regular $4$-ordered hamiltonian
graphs was exhibited. However, the only two $3$-regular $4$-ordered
hamiltonian graphs found were $K_4$ and $K_{3, 3}$.  In this paper we
analyze the class of all $3$-regular graphs with the aim of determining
whether or not there exist other $4$-ordered (and $4$-ordered
hamiltonian)
$3$-regular graphs.  In Section $2$ 
we prove that a $3$-regular $4$-ordered
graph on more than $6$ vertices is not only triangle free \cite{ng}, but it is
also square free and we show that the smallest (by the number of vertices)
$3$-regular triangle and square free graph, namely the Petersen graph, is
$4$-ordered.  We also consider a common family of graphs generalizing the
Petersen graph, and we show that the Petersen graph itself is the only member
of this family that is $4$-ordered.  In Section $3$, we exhibit a
$4$-ordered hamiltonian graph on $14$ vertices,  the Heawood graph, and show
that it is the smallest graph after $K_4$ and $K_{3, 3}$ that is
$3$-regular and $4$-ordered hamiltonian.  In Section $3$ we also exhibit
forbidden subgraphs for $3$-regular and $4$-ordered hamiltonian graphs.
Finally, in Section $4$, we conclude our paper by constructing an infinite family of
$3$-regular $4$-ordered graphs that we call torus-graphs (torus-graphs can be embedded on the torus
without crossing edges). Since $3$-regular graphs have  the lowest
possible degree for $4$-ordered graphs, the construction of torus-graphs
answers the question of whether there are low degree $4$-ordered graphs.

\section{Forbidden subgraphs and the Petersen graph}
\label{s:two}

It is easy to see that no $3$-regular $4$-ordered hamiltonian graph on more
than $4$ vertices contains a triangle  (see \cite{ng}). Also, $3$-regular
$4$-ordered graphs must be triangle free, by an analogous argument.

\begin{thm}
A $3$-regular $4$-ordered graph on more than $6$ vertices  does not contain a square. 
\end{thm}

\begin{proof} Suppose $G$ is a $4$-ordered graph on more than $6$ vertices and
it contains a square. By $4$-orderedness, $G$ is triangle free,
as noted before.  If there exists a square, say with vertices $A$, $B$, $C$, and $D$ (in
order) such that some pair of  edges incident to opposite
vertices of the square $ABCD$ do not share a vertex, then we can show that $G$ is not
$4$-ordered. Indeed, suppose without loss of generality, that the third edge
incident to $A$ is $AE$, that the third edge incident to $C$ is
$CF$, and that $E\neq F$. In this case, there can be no cycle in $G$
containing the vertices \textit{F, E, C, A} in this order because $CF$ and $AE$
cannot be edges in this cycle, which
implies that \textit{AB, BC,
CD, DA} are all edges in this cycle,
which contradicts the existence of a cycle containing vertices \textit{F, E, C, A}
in that order.  
 
We now show that $K_{3, 3}$ is the
only $3$-regular triangle free graph
containing a square such that the edges incident to opposite vertices of any
square it contains do share a vertex. Indeed, suppose $H$ is a
$3$-regular triangle free graph containing a square $ABCD$, such that the edges
incident to opposite vertices of any square it contains do share a vertex. 
  Then, there exists a vertex $E$ and a vertex $F$
such that such that $BE$, $DE$, $AF$, $CF$ are edges.  As $H$ is triangle
free, it follows that $E\neq F$. Consider the square $ADCF$ in $H$, and its
opposite vertices $D$ and $F$.  By assumption, the edges incident to opposite
vertices of any square in $H$  share a vertex, and as the degree of $D$ is
already $3$, it follows that $EF$ is an edge in $H$.  As all of the vertices
$A, B, C, D, E, F$ already have degree $3$ it follows that $H$ is the graph on
these six vertices with edges as described. It is easy to see that $H=K_{3,
3}$.

\end{proof}

\begin{corr}
If $G$ is a $3$-regular $4$-ordered graph on more than $6$ vertices, then 
every vertex has exactly $6$ vertices at distance $2$.  
\end{corr}

\begin{proof} Note that in a $3$-regular graph $G$ any vertex has $2, 3, 4, 5$,
or $6$ vertices at distance $2$.  By Theorem 2.1, in order for graph $G$ on
more than $6$ vertices to be $4$-ordered, it has to be square free. Observe
that if there is a vertex $v$ that has $2, 3, 4$, or $5$ vertices at distance
$2$, then $v$ is a vertex of a  square in $G$. Thus, if $G$ is a $3$-regular
$4$-ordered graph on more than $6$ vertices, then each vertex of $G$ has
exactly $6$ vertices at distance $2$.  
\end{proof}

\begin{figure}   
\begin{center}
\epsfbox{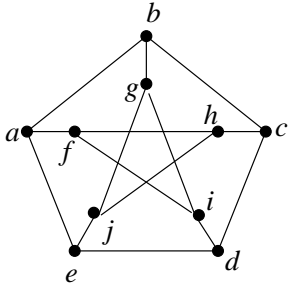}
\caption{The Petersen graph is $4$-ordered (Theorem~\ref{pet4}).}
\label{pet}
\end{center}
\end{figure}

The next two lemmas will be used to prove Theorem~\ref{pet4}, in
which we show that the smallest $3$-regular graph that is triangle and square
free 
is $4$-ordered. This graph is the well-known Petersen graph, (see
Figure~\ref{pet}).

A \textit{walk} is a  sequence of (not necessarily distinct) vertices $x_1,
x_2, \ldots, x_n$ such that $x_i$ is adjacent to $x_{i+1}$ for all $1 \leq i
\leq n-1$.  The length of a walk is the number of edges in the walk.
Following \cite{h}, an \textit{n-route} is a  vertex disjoint walk of length
$n$ with specified initial point. A graph $G$ is \textit{n-transitive}, $n\geq
1$, if it has an $n$-route and if there is always an automorphism of $G$
sending each $n$-route into any other $n$-route.

\begin{lemma} \emph{(\cite{h}, p.175)}
The Petersen graph is $3$-transitive. 
\end{lemma}

The following well-known fact follows from Lemma 2.3:  
 
\begin{lemma} Given any two $5$-cycles in the Petersen graph, there
exists an automorphism that takes one of the $5$-cycles into the other.
\end{lemma}

\begin{thm} \label{pet4}
The Petersen graph is 4-ordered.
\end{thm}

\begin{proof} 
Consider the Petersen graph as shown in Figure~\ref{pet}, and consider $4$
vertices $v_1, v_2, v_3, v_4$ specified in order in the Petersen graph. We
consider three cases depending on how the $4$ vertices are
distributed:
either all  $4$ specified vertices are among vertices \textit{a, b, c, d, e},
or $3$ of them are from \textit{a, b, c, d, e}, or $2$ of them are from
\textit{a, b, c, d, e}. Call the cycle containing the vertices \textit{a, b,
c, d, e} the \textit{outer cycle}, and the cycle containing the vertices
\textit{f, h, j, g, i}  the \textit{inner cycle} of the Petersen graph. 

Consider the case when there are $3$ vertices specified on the outer cycle,
and $1$ vertex specified on the inner cycle. Without loss of generality, the
vertex on the inner cycle can be specified to be the first vertex, $v_1$,
and the $3$ vertices specified on the outer cycle the second, third, and
forth, $v_2, v_3, v_4$. We now show that in this case regardless of 
exactly which $4$ vertices $v_1, v_2, v_3, v_4$ are, there is  a cycle containing
them in this order. 

Let $x_2$ and $x_4$ be the vertices on the inner cycle that are adjacent to
$v_2$ and $v_4$, respectively.  Go from $v_1$ on the inner cycle
\textit{f-h-j-g-i} until $x_2$, without meeting $x_4$. Then go from $x_2$ to
$v_2$, and from $v_2$ go to  $v_3$ and then to $v_4$ on the outer cycle
\textit{a-b-c-d-e}.  From $v_4$ go to $x_4$ and then to $v_1$ without meeting
$x_2$.  This completes the cycle  that contains $v_1, v_2, v_3, v_4$ in this
order. 

Thus, by Lemma 2.4,  if there are exactly $3$ of the $4$ specified vertices on
any $5$-cycle in the Petersen graph, then we have a cycle containing the $4$
vertices in the specified order. 

This observation makes it unnecessary to check the case of all $4$ vertices
being among vertices \textit{a, b, c, d, e}, as in this case there is a
$5$-cycle containing exactly $3$ of the specified vertices. Furthermore, in
the case that  $2$ of the specified vertices are from  \textit{a, b, c, d, e},
and $2$ from \textit{f, g, h, i, j} it suffices to consider the case when
these $4$ vertices are in relative positions as \textit{a, c, i, j} since in
all other cases there is a $5$-cycle containing $3$ of the specified vertices.
For these remaining cases, one can easily find a cycle containing the vertices
no matter how we specify their order.

\end{proof}

We now  consider a common family of graphs generalizing the Petersen graph and
show that the Petersen graph  itself is the only member of this family that is
$4$-ordered, although they are all $3$-ordered hamiltonian with the exception of 
the Petersen graph (\cite{hol}, p.136). 

A \textit{star graph} $S_{n, k}$, where $n, k$ are positive integers, is a
graph  on  vertices $A_1, A_2, \ldots, A_n$ such that $A_i$ and $A_j$ are
adjacent if the indices $i$ and $j$  differ by $k$ modulo $n$.  The
\textit{generalized Petersen graph} $P_{n, \lfloor \frac{n-1}{2}\rfloor}$,
$n\geq 5$, is a graph consisting of an cycle of length  $n$ on the vertices
$B_1$, $B_2$, \ldots, $B_n$ (in order)  and a star
graph $S_{n, \lfloor \frac{n-1}{2}\rfloor}$ on vertices $A_1, A_2, \ldots,
A_n$, such that
$A_i$ and $B_i$ are adjacent for all $i=1, 2, \ldots, n$.  One can imagine the vertices  $B_1,B_2,\ldots, B_n$ to be drawn on an
outer cycle and the vertices $A_1, A_2, \ldots, A_n$ on an inner cycle.  Note
that $P_{5,2}$ is the standard Petersen graph.

\begin{prop} \label{genpet4}
The generalized Petersen graph $P_{n, \lfloor \frac{n-1}{2}\rfloor}$ is not $4$-ordered for $n>5$. 
\end{prop}

\begin{figure}   
\begin{center}
\epsfbox{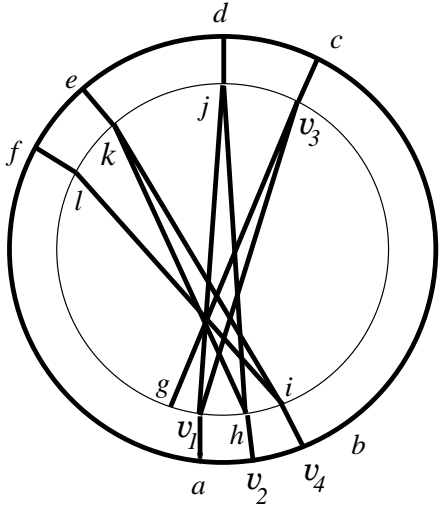}
\caption{The generalized Petersen graphs are not $4$-ordered (Proposition~\ref{genpet4}).}
\label{genpet}
\end{center}
\end{figure}

\begin{proof} If $n=6$, $P_{6, 2}$ is $3$-regular and contains a triangle,
thus it is not $4$-ordered.  In   Figure~\ref{genpet}  we depicted some edges
of $P_{n, \lfloor \frac{n-1}{2}\rfloor}$ in the case $n>6$.  The thick lines
denote edges in the graph, and the thin line denotes the circle upon which the
star graph is drawn. Note that if $n=7$, then  $b=c$ and $g=l$ in
Figure~\ref{genpet}; however, this has no effect on the following arguments.  

In Figure~\ref{genpet} we also marked vertices \textit{1, 2, 3, 4} and we now
show  that there is no cycle containing these vertices in this specified
order. Indeed, suppose that there is such a cycle $C$.  As \textit{13} and
\textit{24} are edges that cannot occur in $C$, it follows that $C$ must
contain edges \textit{j1, 1a, a2, 2h, i4, 4b, g3,} and \textit{3c}; and since
\textit{1-a-2}  is  a path in $C$  this also shows that the edges \textit{j1,
1a, a2}, and \textit{2h} are oriented from their first vertex to the second
(i.e., $j1$ is
oriented from $j$ to $1$, etc.). As \textit{hk} and \textit{hj} are the
remaining  edges from \textit{h}, and  \textit{j1} is an edge in $C$, it follows that $C$ contains \textit{hk} with this orientation.
Also, as \textit{il} and \textit{ik} are the remaining  edges from \textit{i}, and
\textit{k} has been visited when we were going from \textit{2} to \textit{3},
it follows that $C$ contains the edge \textit{il}. As  \textit{kh, ki}, and \textit{ke}
are the edges from \textit{k}, and edge $hk$ has been used, and $ki$ cannot be
used as it would create a path in $C$ directly from \textit{2} to \textit{4},
it follows that $C$ contains the edge \textit{ke} with this orientation.
Because edge $hj$ cannot be in $C$, the edge $dj$ must be in $C$, with this
orientation.  If the edge $cd$ were in $C$, then there would be a direct path
in $C$ from \textit{3} to \textit{1}, which is impossible; and
thus $ed$ must be an edge in $C$.  But then there is a closed cycle from
\textit{1} to \textit{2} and back to \textit{1} in $C$, which is also
impossible, giving us a contradiction.
\end{proof}

\section{The Heawood Graph} 
 
\begin{figure}    
\begin{center}
\epsfbox{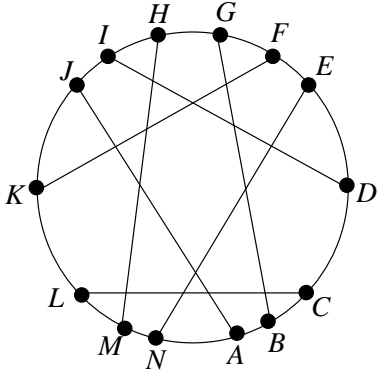}
\caption{The Heawood graph is $4$-ordered hamiltonian (Theorem~\ref{hea4}).}
\label{hea}
\end{center}
\end{figure}

In this section we show  that the Heawood graph,  the smallest $3$-regular
graph that is triangle, square, and pentagon free (\cite{chart}, p. 42), is
$4$-ordered hamiltonian (Figure~\ref{hea}).

\begin{lemma} \label{l:h4trans}
\emph{(\cite{h}, p.174)} The Heawood graph is 4-transitive.
\end{lemma}

\begin{corr} \label{diam3} The diameter of the  Heawood graph is 3.
Furthermore, if two vertices $X$ and $Y$ of the  Heawood graph are at distance
3, then there are two disjoint paths of length 3 between $X$ and $Y$.
\end{corr}

\begin{proof}  

Take any vertex $v$ in the Heawood graph. Observe that in order to prove that
the diameter is $3$, it suffices to check that all vertices are at distance
less than $4$ from $v$ and there is a vertex at distance $3$ from $v$, which
is easily done.  In order to show that if two vertices $X$ and $Y$ of the
Heawood graph are at distance 3, then there are two disjoint paths of length
$3$ between $X$ and $Y$, it suffices to find two vertices $X$ and $Y$ at
distance $3$ with two disjoint paths between them, say $P_1$ and $P_2$, which
is easily done. Then, by applying Lemma 3.1 the claim follows for any two
vertices $X'$ and $Y'$ at distance $3$. 
\end{proof}

\begin{thm} \label{hea4}
The Heawood graph is 4-ordered hamiltonian.
\end{thm}
 
\begin{proof}
We will consider cases depending on the
distances between the $4$ specified vertices $v_1, v_2,
v_3, v_4$. We denote the distance between vertices $X$ and $Y$ by $d(X, Y)$. 

\textit{Case 1.} Some $2$ of the $4$ specified vertices are at distance $1$.
Without loss of generality these are either $v_1$ and $v_2$ or $v_1$ and
$v_3$.

\textit{Case 1.1.1.} $d(v_1, v_2)=1$, $d(v_2, v_3)=1$. By
Lemma~\ref{l:h4trans}, we can suppose without loss of generality that $v_1,
v_2, v_3$ are respectively \textit{A, B, C}. In this case wherever $v_4$ is,
it is clear that \textit{A-B-C-D-E-F-G-H-I-J-K-L-M-N-A} is the desired
hamiltonian cycle. 

\textit{Case 1.1.2.} $d(v_1, v_2)=1$, $d(v_2, v_3)=2$, and there is a length
$2$ path from $v_2$ to $v_3$ not containing $v_1$.  By Lemma~\ref{l:h4trans},
we can suppose without loss of generality that $v_1, v_2, v_3$ are
respectively \textit{A, B, D}. In the case $v_4\neq C$, it is clear that
\textit{A-B-C-D-E-F-G-H-I-J-K-L-M-N-A} is the desired hamiltonian cycle. In
case $v_4 = C$, the desired hamiltonian cycle is
\textit{A-B-G-H-I-J-K-F-E-D-C-L-M-N-A}. 

\textit{Case 1.1.3.} $d(v_1, v_2)=1$, $d(v_2, v_3)=2$, and there is no length
$2$ path from $v_2$ to $v_3$ not containing $v_1$.  We can suppose  that $v_1,
v_2, v_3$ are respectively \textit{A, B, N}.  In the case $v_4 \neq C, M, L$, the
desired hamiltonian cycle is \textit{A-B-C-L-M-N-E-D-I-H-G-F-K-J-A}.  In the case
$v_4=C$ or $v_4=M$ or $v_4 = L$, the desired hamiltonian cycle is
\textit{A-B-G-F-E-N-M-H-I-D-C-L-K-J-A}. 
 
\textit{Case 1.1.4.} $d(v_1, v_2)=1$, $d(v_2, v_3)=3$. By Corollary 3.2 there
is a  length $3$ path from $v_2$ to $v_3$ not containing $v_1$. Thus, by
Lemma~\ref{l:h4trans} we can suppose  that $v_1, v_2, v_3$ are respectively
\textit{A, B, E}.

In the case $v_4\neq F, G, H, I, J, K$, the desired hamiltonian cycle is
\textit{A-B-G-H-I-J-K-F-E-D-C-L-M-N-A}.  Clearly, in the case $v_4$ is \textit{F,
G, H, I, J}, or $K$ the desired hamiltonian cycle is
\textit{A-B-C-D-E-F-G-H-I-J-K-L-M-N-A}. 
 
\textit{Case 1.2.1.} $d(v_1, v_3)=1$, $d(v_3, v_4)=1$. We can suppose by Lemma
3.1 that $v_1, v_3, v_4$ are respectively \textit{A, B, C}. If $v_2=D, E, F,
G, H, I, J$, or $K$, the desired hamiltonian cycle is
\textit{A-J-K-F-E-D-I-H-G-B-C-L-M-N-A}. On the other hand, if $v_2=L$ or
$v_2=M$, according to Lemma~\ref{l:h4trans} this is the same as if $v_1, v_2,
v_3, v_4$ are respectively \textit{A, D, B, C} or \textit{A, E, B,
C}; thus
this is also covered by the previous case.  Finally, if $v_2=N$, then the
desired hamiltonian cycle is \textit{A-N-E-F-K-L-M-H-G-B-C-D-I-J-A}.

\textit{Case 1.2.2.} $d(v_1, v_3)=1$, $d(v_3, v_4)=2$, and there is a length
$2$ path from $v_3$ to $v_4$ not containing $v_1$.  We can suppose that $v_1,
v_3, v_4$ are respectively \textit{A, B, D}. In the case $v_2=C$, the desired
hamiltonian cycle is \textit{A-J-K-L-C-B-G-F-E-D-I-H-M-N-A}. In the case $v_2\neq
C, I, J$, the desired hamiltonian cycle is
\textit{A-N-E-F-K-L-M-H-G-B-C-D-I-J-A}.  In the case $v_2=I$ or $v_2=J$ the
desired hamiltonian cycle is \textit{A-J-I-H-M-L-K-F-G-B-C-D-E-N-A}.

\textit{Case 1.2.3.} $d(v_1, v_3)=1$, $d(v_3,v_4)=2$, and  every length $2$
 path from $v_2$ to $v_3$ contains $v_1$.  We can suppose  that $v_1, v_3,
 v_4$ are respectively \textit{A, B, N}.  In the case $v_2\neq E, F, G$, the
 desired hamiltonian cycle is \textit{A-J-K-L-M-H-I-D-C-B-G-F-E-N-A}.  In the case
 $v_2=E$ or $v_2=F$ or $v_2=G$, the desired hamiltonian cycle is
 \textit{A-J-K-F-E-D-I-H-G-B-C-L-M-N-A}. 
 
\textit{Case 1.2.4.} $d(v_1, v_3)=1$, $d(v_3,v_4)=3$. By
 Corollary~\ref{diam3} there is a path from $v_3$ to $v_4$ not containing
 $v_1$. Thus, we can suppose  that $v_1, v_3, v_4$ are respectively \textit{A,
 B, E}. Applying Lemma~\ref{l:h4trans} to $3$-routes \textit{B-A-N-E} and
 \textit{A-B-C-D},  we can identify this case with $v_3, v_1, v_4$ being
 \textit{A, B, D} respectively. If $v_2=C, L, M$, or $N$, the desired
 hamiltonian cycle is \textit{B-C-L-M-N-A-J-K-F-E-D-I-H-G-B}.  If $v_2=F, G,
 J$, or $K$ the desired hamiltonian cycle is
 \textit{B-G-F-K-J-A-N-E-D-I-H-M-L-C-B}.  The remaining cases that we have to
 consider are when $v_2=E, H$, or $I$. As \textit{A-B-C-D-E} and
 \textit{A-B-C-D-I} are both $4$-routes, by Lemma~\ref{l:h4trans} it suffices
 to consider only the cases when $v_2=H$ or $v_2=I$. In this case the desired
 hamiltonian cycle is \textit{B-G-H-I-J-A-N-M-L-K-F-E-D-C-B}.

\textit{Case 2.} No $2$ of the $4$ specified vertices are at distance $1$. As
the diameter is $3$, the possible  distances are $2$ and $3$. 

\textit{Case 2.1.} $d(v_1, v_2)=d(v_1, v_3)=d(v_1, v_4)=2$. Without loss of
generality we can assume that either $v_4, v_1, v_2$ are \textit{N, B, D} or
$v_4, v_1, v_2$ are \textit{L, B, D}.  

In the case $v_4, v_1, v_2$ are \textit{N, B, D}, as no two of the four specified
vertices are at distance $1$, $v_3\neq A, C$; and thus the
desired hamiltonian cycle is \textit{B-C-D-E-F-G-H-I-J-K-L-M-N-A-B}. 

In the case $v_4, v_1, v_2$ are \textit{L, B, D}, as no two of the four specified
vertices are at distance $1$, $v_3\neq A, C, M$; and thus if $v_3\neq N$ the desired
hamiltonian cycle is \textit{B-C-D-E-F-G-H-I-J-K-L-M-N-A-B}.  If $v_3=N$, then
the desired hamiltonian cycle is \textit{B-G-F-E-D-I-H-M-N-A-J-K-L-C-B}.

\textit{Case 2.2.} Some $2$ vertices are at distance $2$ from $v_1$, and $1$
is at distance $3$.  The case when $v_2$ and $v_4$ are at distance $2$ from
$v_1$ can be solved analogously to Case 2.1.  Consider the case $d(v_1,
v_2)=d(v_1, v_3)=2$ and $d(v_1, v_4)=3$. Without loss of generality, either
$v_3, v_1, v_2$ are \textit{N, B, D} or $v_3, v_1, v_2$ are \textit{L, B, D}.
In the case $v_3, v_1, v_2$ are \textit{N, B, D}, the only
possibility for $v_4$ is
$K$, and the desired hamiltonian cycle is
\textit{B-C-D-E-N-M-L-K-F-G-H-I-J-A-B}.  In the case $v_3, v_1, v_2$ are
\textit{L, B, D}, then all the points that are distance $3$ from $v_1$ would
be at distance $1$ from some of the vertices $v_2$ or $v_3$, contradicting our
assumption. 

\textit{Case 2.3.} Some $2$ vertices are at distance $3$ from $v_1$, and the
remaining vertex is at distance $2$.  We can suppose without loss of
generality that either $d(v_1, v_2)=2$, or $d(v_1, v_3)=2$.  

If $d(v_1, v_2)=2$, suppose $v_1, v_2$ are \textit{B, D}. As the $4$ vertices
at distance $3$ from $B$ are \textit{E, I, K, M} and $d(v_2, E)=1$, $d(v_2,
I)=1$, we have that $\{v_3, v_4\}=\{K, M\}$. If $v_4=K$ and $v_3=M$, then the
desired hamiltonian cycle is \textit{B-C-D-E-N-M-L-K-F-G-H-I-J-A-B}, whereas
if $v_3=K$ and $v_4=M$, then the desired hamiltonian cycle is
	\textit{B-C-D-E-F-G-H-I-J-K-L-M-N-A-B}. 

If $d(v_1, v_3)=2$, suppose $v_1, v_3$ are \textit{B, D}. By an analogous
argument as above $\{v_2, v_4\}=\{K, M\}$. If $v_2=K$ and $v_4=M$, then the
desired hamiltonian cycle is \textit{B-A-J-K-F-G-H-I-D-E-N-M-L-C-B}, whereas
if $v_4=K$ and $v_2=M$, then the desired hamiltonian cycle is
	\textit{B-A-N-M-H-G-F-E-D-I-J-K-L-C-B}. 

\textit{Case 2.4.} Suppose that all $3$ points are at distance $3$ from $v_1$.
Without loss of generality, let $v_1=B$. The distance $3$ vertices from $v_1$
are \textit{E, I, K, M} and regardless of which $3$ of these we choose for
$v_2, v_3, v_4$, there will be two at distance 2.  However, as the choice of
$v_1$ was without loss of generality in all previous cases,  this case cannot
occur.

\end{proof}

\begin{thm} The Heawood graph is the graph on the fewest vertices, after $K_4$ and
$K_{3, 3}$, that is $3$-regular $4$-ordered hamiltonian.
\end{thm}

\begin{proof} The only $3$-regular $4$-ordered hamiltonian graphs on fewer than
$7$ vertices are $K_4$ and $K_{3, 3}$, \cite{ng}. By \cite{ng} and Theorem
2.1, all $3$-regular $4$-ordered hamiltonian graphs on more than $6$ vertices
are both triangle and square free. The smallest such $3$-regular graph is the
Petersen graph on $10$ vertices, and it is not hard to see that it is the only
$3$-regular triangle- and square-free graph on $10$ vertices.  As it is not
hamiltonian, it is not $4$-ordered hamiltonian either. Since a $3$-regular
graph has even number of vertices, the only  possibility for a $3$-regular $4$-ordered hamiltonian
graph on fewer than $14$ vertices is a $3$-regular graph on $12$ vertices that
is triangle and square free.

\begin{figure}
\begin{center}
\epsfbox{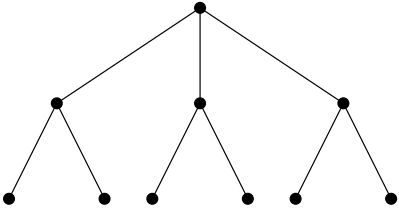}
\caption{Each vertex has six vertices at distance two.}
\label{6dist2}
\end{center}
\end{figure}

As square-freeness of a graph $G$ is equivalent to each vertex
of $G$ having $6$ vertices at distance $2$, it follows 
given any triangle- and square-free graph $G$ on at least $10$ vertices 
and an arbitrary vertex
$v$ of $G$, the induced subgraph of $G$ on $v$ and the vertices that are at
distance $1$ and at distance $2$ from $v$ contains the graph in Figure~\ref{6dist2} as a
subgraph.   Since the graph in Figure~\ref{6dist2} has $10$ vertices, only two
vertices have to be added to obtain a $3$-regular triangle- and square-free graph on $12$ vertices. It is not hard
to see that there are only two non-isomorphic $3$-regular graphs on $12$ vertices
that are triangle free and square free (Figure~\ref{12vert}).

\begin{figure}
\begin{center}
\epsfbox{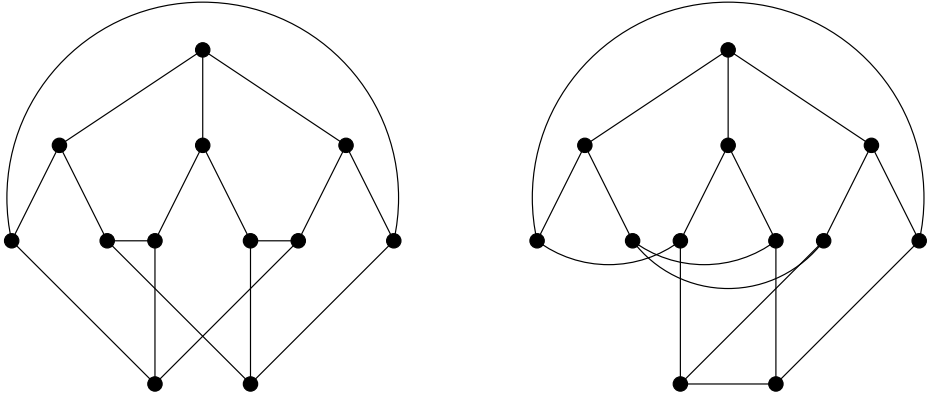}
\caption{The two triangle- and square-free $3$-regular graphs on $12$ vertices.}
\label{12vert}
\end{center}
\end{figure}

We now show that neither of the two triangle- and square-free $3$-regular graphs
on $12$ vertices are $4$-ordered hamiltonian. 

\begin{figure}
\begin{center}
\epsfbox{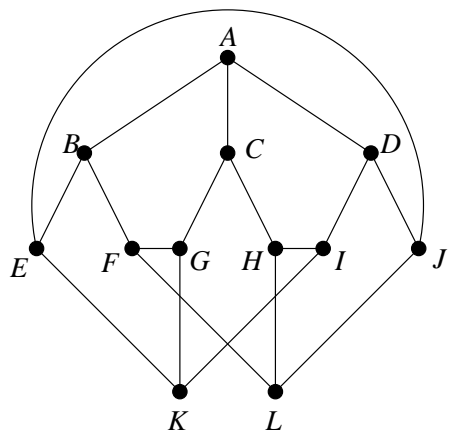}
\caption{There is no hamiltonian cycle containing vertices $A, E, D, J$ in this order.}
\label{elso}
\end{center}
\end{figure}

Consider Figure~\ref{elso}. We show that there is no hamiltonian cycle
containing vertices $A, E, D, J$ in this order. Suppose for sake of
contradiction that there is a
hamiltonian cycle $\mathcal{C}$ containing vertices $A, E, D, J$ in this
order. Then edges $AD$ and $EJ$ are not in the cycle, thus edges $CA, AB, BE,
EK, ID, DJ, JL$ are in the cycle, and it is not hard to see
that they are
oriented in this way (i.e.,  we are going from $C$ to $A$, from $A$ to $B$,
etc.). Note that the only way to get from $E$ to $D$ without meeting $C$ and
$L$ (vertices that are on $\mathcal{C}$ after meeting all of $A, E, D, J$) is
to take the edge $KI$ from $K$ to $I$. However, there is no path from $L$ to
$C$ exactly including vertices $G$ and $H$, thus there is no hamiltonian cycle
$\mathcal{C}$ containing vertices $A, E, D, J$ in this order.

\begin{figure}
\begin{center}
\epsfbox{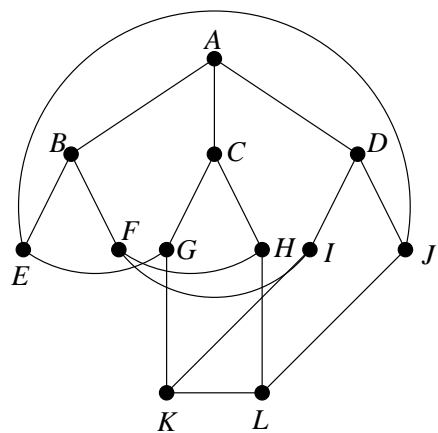}
\caption{There is no hamiltonian cycle containing vertices $A, E, C, G$ in this order.}
\label{masodik}
\end{center}
\end{figure}

Consider Figure~\ref{masodik}. We show that there is no hamiltonian cycle
containing vertices $A, E, C, G$ in this order. Suppose for sake of
contradiction that there is a
hamiltonian cycle $\mathcal{C}$ containing vertices $A, E, C, G$ in this
order.  Then edges $AC$ and $EG$ are not in the cycle, thus edges $DA, AB,
BE, EJ, HC, CG, GK$ are in the cycle, and it is not hard to see that they
are oriented in this way (i.e.,  we are going from $D$ to $A$, from $A$ to
$B$, etc.). Since edge $JD$ cannot be in $\mathcal{C}$ as it would mean that
$\mathcal{C}$ does not contain $C$ and $G$, $\mathcal{C}$ contains edge $ID$,
with this orientation. If $KI$ were not in $\mathcal{C}$, this would imply $FI$
and $HF$ were in $\mathcal C$, which is impossible given that $HC$ is oriented from $H$ to $C$.
Thus, edge $KI$ is in $\mathcal{C}$. As $DJ$ is not in $\mathcal{C}$, it
follows that $JL$ is in $\mathcal C$. Now we
have to get from $L$ first to $C$ and then to $G$, and so
it follows that edge $LH$ is in $\mathcal{C}$. However, in this case
$\mathcal{C}$ does not contain $F$, therefore $\mathcal{C}$ is not a
hamiltonian cycle, leading to the desired contradiction.

\end{proof}

\begin{figure}
\begin{center}
\epsfbox{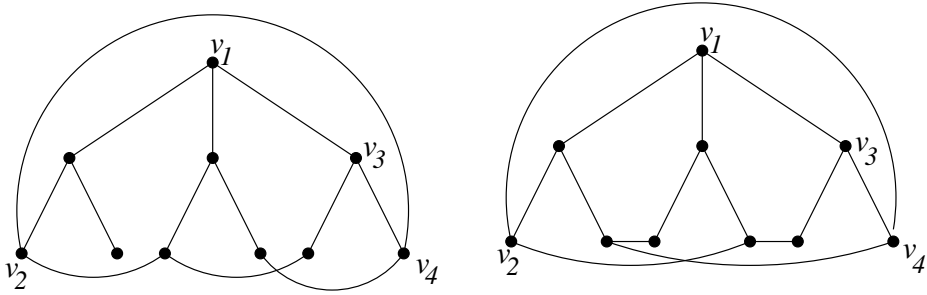}
\caption{There is no $4$-ordered hamiltonian graph $G$ on more than $10$ vertices that contains 
these graphs as subgraphs.}
\label{forb}
\end{center}
\end{figure}

\begin{prop} A $4$-ordered hamiltonian graph $G$ on more than $10$ vertices
contains neither of the graphs in Figure~\ref{forb} as subgraphs. 
\end{prop}

\begin{proof} The proof 
follows the arguments in the proof of Theorem 3.4, where $v_1, v_2,
v_3, v_4$ are the four vertices for which there is no hamiltonian cycle
containing them in this order. See Figure~\ref{forb}.
  \end{proof}

\begin{figure}
\begin{center}
\epsfbox{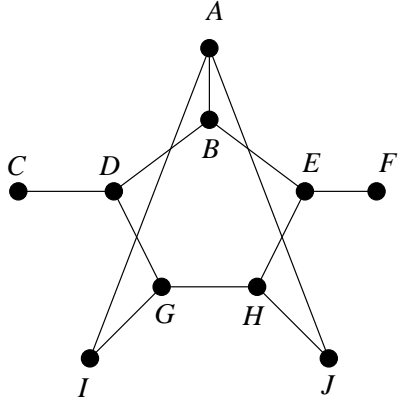}
\caption{There is no $4$-ordered hamiltonian graph $G$ that contains this graph as 
a subgraph.}
\label{last}
\end{center}
\end{figure}

\begin{prop}
There is no $4$-ordered hamiltonian graph $G$ that contains the graph in 
 Figure~\ref{last} as a subgraph.
\end{prop}

\begin{proof}
Suppose the graph in Figure~\ref{last}
 is a subgraph of a $4$-ordered hamiltonian graph $G$. Then there exists a hamiltonian 
cycle $\mathcal{C}$ in $G$ containing vertices \textit{D, E, G, H} in this order. Since edges
 \textit{DG} and \textit{EH} cannot be in $\mathcal{C}$, edges \textit{CD, DB, BE, EF, IG, GH} 
 and \textit{HJ} are in $\mathcal{C}$. As $\mathcal{C}$ is a hamiltonian cycle, and since edge 
\textit{BA} cannot be in $\mathcal{C}$ as \textit{DB} and \textit{BE} are, it follows that edges 
\textit{IA} and \textit{AJ} must be in $\mathcal{C}$. However, edges \textit{AI, IG, GH, HJ}
and \textit{JA} form a cycle, which contradicts the existence of $\mathcal{C}$. 
 
\end{proof}

\section{A Family of 4-ordered 3-regular graphs}

\begin{figure}    
\begin{center}
\epsfbox{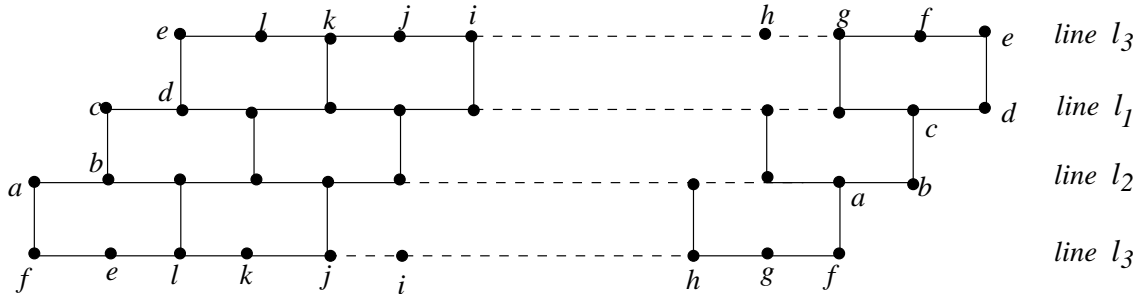}
\caption{The torus-graphs.
 Note that the vertices with the same names are
identified, and consequently  $l_3$ is both the bottom and the top line.}
\label{torus}
\end{center}
\end{figure}

In this section we introduce an infinite family of $3$-regular $4$-ordered
graphs that answers the question of the existence of an infinite family of low
degree $4$-ordered graphs.  These graphs, which we will refer to as {\em
torus-graphs}, are similar to the Heawood graph
in that they have girth $6$, and just as the Heawood graph can be embedded on the torus, they
too can be embedded  in the torus.
 A general
torus-graph is shown in Figure~\ref{torus}.  The
dashed lined signify that there are more hexagons in each of the three
rows,
and we will always assume we have sufficiently many hexagons in each row.
The labelings \textit{a, b, c, ..., k, l} show which vertices are identified.  

\begin{thm} The torus-graphs presented in Figure~\ref{torus} are $3$-regular
$4$-ordered graphs for long enough rows.  \end{thm}

\begin{proof}

We  consider the graphs as drawn in the plane, and we shall refer to left and
right movements, as well as up and down movements, when
referring to moves on
the horizontal and vertical lines respectively. Furthermore, we refer to the
horizontal lines in Figure~\ref{torus} as lines $l_1$, $l_2$, $l_3$, where the
top and the bottom lines are $l_3$, and $l_1$ and $l_2$ are,
respectively,  the second and
third horizontal lines from the top.

We analyze all possible ways of specifying the $4$ points, $v_1, v_2, v_3,
 v_4$, and in each case we exhibit a cycle containing them in the specified
 order. 

We consider cases depending on how the $4$ points are
distributed among the lines. 

\textit{Case 1.} All $4$ points are on the same line.  
Without loss of generality, say the line is $l_2$. 

\textit{Case 1.1.} If beginning from $v_1$ the $4$ points are
in  the order $v_1$, $v_2$, $v_3$, $v_4$
along the line, either going to the right on the
line or to the left starting at $v_1$, then   just trace along $l_2$
to obtain the desired cycle. 

\textit{Case 1.2.} If the $4$ points are not in  order $v_1, v_2, v_3, v_4$,
 we can suppose without loss of generality that they are in order $v_1, v_2,
 v_4, v_3$ when going to the right from $v_1$.  In this case we consider two
 subcases.

\textit{Case 1.2.1.} $v_1$ and $v_3$ are not adjacent. 
 Then, do the following: 

\begin{itemize}

\item go right from $v_1$ to $v_2$ on line $l_2$

 \item go down or up from $v_2$ (depending whether the vertical edge starting
at $v_2$ is down or up) to the nearest line; suppose without loss of
generality that the edge we took was down
 
 \item go right until reaching either the vertical edge that
 has $v_3$ as its end  or
the vertical edge with an end that is one step
to the right of $v_3$ , and traverse this vertical
edge, returning to the line $l_2$
  
 \item go left, hitting $v_3$ and $v_4$ and go until the vertex   next
to $v_2$, and from this vertex take the edge up
 
 \item if the vertical edge from $v_1$ is up, continue on this line ($l_1$)
until hitting that edge and take it to $v_1$; if the vertical edge from $v_1$
is down, go one step left, go up, and then continue left along $l_3$ 
 until hitting the vertical
edge of $v_1$ and take it to $v_1$
 
 \end{itemize}

\textit{Case 1.2.2.} $v_1$ and $v_3$ are adjacent. We can suppose without loss
of generality that the vertical edge of $v_1$ is down, and
consequently that the
vertical edge of $v_3$ is up.
 
 If the vertical edge of $v_2$ is down, 
\begin{itemize}
 
\item go right from $v_1$ to $v_2$

\item take the vertical edge from $v_2$ down to $l_3$

\item go right on $l_3$ until the vertex that is next to the vertical edge from $v_1$

\item go down to $l_1$

\item take a step left, and then  step down to $l_2$, hitting $v_3$

\item go left, hitting $v_4$, until the vertex   next to $v_2$, and from
this vertex take the edge up
 
\item go one step left, go up, and then go left until hitting the vertical
edge of $v_1$ and take it to $v_1$
  
\end{itemize}

 If the vertical edge of $v_2$ is up, 
\begin{itemize}
 
\item go right from $v_1$ to $v_2$

\item take the vertical edge from $v_2$ up to $l_1$

\item go right until hitting the vertical edge from $v_3$

\item go down to $l_1$ hitting $v_3$

\item go left, hitting $v_4$, until the vertex that is next to $v_2$, and from
this vertex take the edge down
 
\item go left until hitting the vertical edge of $v_1$ and take it to $v_1$
  
\end{itemize}

\textit{Case 2.} There are $3$ points on one line (say on line $l_1$) and $1$
on another. Without loss of generality we can suppose that the 3 points on
line $l_1$ are $v_1, v_2, v_3$, and that this is the order as we go from $v_1$
to the right. Also, we can suppose without loss of generality that $v_4$ is on
line $l_2$.

\textit{Case 2.1.} Both $v_1$ and $v_3$ have vertical edges down.

\begin{itemize}

\item go from $v_1$ to the right to $v_2$ and then to $v_3$

\item go down, go in the direction so as to hit $v_4$ and then the vertical edge of $v_1$

\item go up to $v_1$
 
\end{itemize}

\textit{Case 2.2.} Both $v_1$ and $v_3$ have vertical edges up. Analysis can
be done similarly as before.

\textit{Case 2.3.} If the vertical edge from $v_1$ is up, from $v_3$ is down,
$v_1$ and $v_3$ are non-adjacent,
\begin{itemize}
 
\item go from $v_1$ to the right to $v_2$ and then to $v_3$

\item go down, go in the direction so as to hit $v_4$ and then the vertical
edge to the left of $v_1$

\item go up to $l_1$ and go one step right to $v_1$ 
 
\end{itemize}

If the vertical edge from $v_1$ is up, from $v_3$ is down, $v_1$ and $v_3$ are
 adjacent,

\begin{itemize}
\item go from $v_1$ to the right to $v_2$ and then to $v_3$

\item go down, go in either  direction so as to hit $v_4$ and then take the first
vertical edge down to $l_3$

\item go on $l_3$ so as to hit the vertical edge of $v_1$, and then take it
 and go to $v_1$
 
\end{itemize}

\textit{Case 2.4}  If the vertical edge from $v_1$ is down, and the vertical edge 
 from $v_3$ is up,
$v_1$ and $v_3$ are non-adjacent, or if the vertical edge from $v_1$ is down,
and the vertical edge from $v_3$ is up, $v_1$ and $v_3$ are adjacent, an analysis similar as before 
  gives the desired cycle.

\textit{Case 3.} There is a line with exactly $2$ vertices on it. Without loss
 of generality these are either $v_1$ and $v_2$ or $v_1$ and $v_3$.

\textit{Case 3.1.} $v_1$ and $v_2$ are on line $l_1$.  The construction of a
cycle is easy in this case, considering two cases depending on whether $v_3$
and $v_4$ are on the same line (without loss of generality on $l_2$) or on 
two different lines.

\textit{Case 3.2.} $v_1$ and $v_3$ are on the same line $l_1$ and they are
non-adjacent.

\textit{Case 3.2.1} $v_2$ and $v_4$ are on the same line, without loss of
generality on $l_2$ 
 
If $v_2$ and $v_4$ are non-adjacent, 

\begin{itemize}

\item go from $v_1$ down if there is a vertical edge, or if not go to one of
its neighbors, so that when taking the edge down from the neighbor, the path
is shorter to $v_2$ (without hitting $v_4$)

\item go to $v_2$, and then go up to $l_1$  immediately or after the next step

\item go along $l_1$ until hitting $v_3$ (without hitting $v_1$) and go down either
from $v_3$ or from the next vertex

\item go until hitting $v_4$ (without hitting $v_2$) and then go down  to $l_3$ either
 from $v_4$ or from the next vertex

\item go along $l_3$ until hitting the vertical from $v_1$, or the vertical from one of
the neighbors of $v_1$ that was not yet used

\item take the vertical edge to $l_1$ and go to $v_1$

\end{itemize}

If $v_2$ and $v_4$ are adjacent, a similar analysis can be performed, also
taking into account whether the vertical edges from $v_1, v_2, v_3$, and
$v_4$ are up or down.

\textit{Case 3.2.2} $v_2$ and $v_4$ are on different lines, without loss of
generality $v_2$ on $l_2$ and $v_4$ on $l_3$. 

\begin{itemize}

\item either go down from $v_1$, or take one step left and go down to line $l_2$

\item go in one of the directions along $l_2$ until hitting $v_2$ so
that either it is
possible to go up from $v_2$, or it is possible to continue one more step and
go up to $l_1$

\item go until hitting $v_3$ (in the direction along $l_1$ so that $v_1$
is not hit again), and
go up either at $v_3$ or continue one more step and go up to $l_3$

\item go in the direction along $l_3$ so as to hit
$v_4$ and then hit either the vertical edge leading to $v_1$ or the vertical
edge leading to one of the neighbors of $v_1$; take this vertical edge and go to $v_1$ 

\end{itemize}

\textit{Case 3.3.} $v_1$ and $v_3$ are on the same line $l_1$ and they are
adjacent.  

\textit{Case 3.3.1} $v_2$ and $v_4$ are on the same line. 

If $v_2$ and $v_4$ are non-adjacent, the analysis can be carried out in a way
similar to  Case 3.2.1. 

If $v_2$ and $v_4$ are adjacent, the analysis can be carried out depending on
 whether the vertical edges from $v_1, v_2, v_3$, and
$v_4$ are up or down, 
taking into account the relative position of $v_1, v_3$ and $v_2, v_4$.

\textit{3.3.2.} $v_2$ and $v_4$ are on different lines.
Without loss of
generality $v_2$ is on $l_2$ and $v_4$ on line $l_3$

 If the vertical edge from $v_1$ is down and the vertical edge from $v_3$ is up, 
\begin{itemize}

\item go from $v_1$ down to $l_2$

\item go to $v_2$ and take the vertical edge from $v_2$ or from a neighbor up to $l_1$

\item go to $v_3$ and take the vertical edge up to $l_3$

\item go in one of the directions along $l_3$ so  so as to hit first
$v_4$ and then the vertical edge that is one step to the left of $v_1$

\item take the vertical edge to $l_1$ and take one step right, ending
the cycle at $v_1$
 
\end{itemize}

If the vertical edge from $v_1$ is up and the vertical edge from $v_3$ is down,
the analysis can be carried out in a similar fashion. 
  
\end{proof}

\section{Acknowledgments} 
This research was performed at the University of Minnesota Duluth under the supervision of 
Professor Joseph Gallian. The author would like to thank Professor Gallian for his support 
and encouragement as well as Denis Chebikin, 
Philip Matchett, David Moulton and Melanie Wood for many useful
 suggestions. Financial support was provided by the Massachusetts Institute of Technology and 
 the University of Minnesota Duluth.

\end{document}